\documentclass{amsart}
\usepackage{pdfpages, 
latexsym,amssymb}
\def\proof{\medskip\noindent{\sc Proof. }}
\def\EOP{\hfill$\Box$}

\def\reals{{\mathbb R}}

\def\field{\hbox{\rm F\kern-.5em F}}
\def\proj{\hbox{\rm P\kern-.5em P}}
 
\newtheorem{thm}{Theorem}[section]
\newtheorem{lem}[thm]{Lemma}
\newtheorem{prop}[thm]{Proposition}
\newtheorem{defn}[thm]{Definition}
\newtheorem{cor}[thm]{Corollary}
\newtheorem{qn}{Question}[section]
\newtheorem{examp}[thm]{Example}

\newtheorem{rk}[thm]{Remark}

\begin{document}
\title{On optimizing discrete Morse functions}
\author{Patricia Hersh}
\curraddr{
        Patricia Hersh,
        Department of Mathematics,
        Indiana University,
        Rawles Hall,
        Bloomington, IN 47405}        
\email{phersh@indiana.edu}
\thanks{The author was supported by an NSF Postdoctoral Research
Fellowship during part of this work.}
\maketitle
 
\begin{abstract}
In 1998, Forman introduced discrete Morse theory as a tool for 
studying CW complexes by producing
smaller, simpler-to-understand complexes of critical cells with the same
homotopy types as the original complexes.
This paper addresses two questions: (1) under what conditions may several 
gradient paths in a discrete Morse function simultaneously be reversed to 
cancel several pairs of critical cells, to further collapse the complex,
and (2) which gradient paths are individually reversible in
lexicographic discrete Morse functions on poset order complexes.  The
latter follows from a correspondence between gradient paths and
lexicographically first reduced expressions for permutations.  
As an application, a new partial order on the symmetric 
group recently introduced by Remmel is proven to be Cohen-Macaulay.
\end{abstract}

\section{Introduction}\label{intro}

Forman introduced discrete Morse theory in [10] as a tool for 
studying the homotopy type and homology groups of 
finite CW-complexes.  In joint work with Eric Babson, 
we introduced a way of constructing ``lexicographic
discrete Morse functions'' for the order complex of any finite poset with
$\hat{0} $ and $\hat{1}$ in [2].
Lexicographic discrete Morse functions have relatively few critical 
cells: when one builds the complex by sequentially attaching
facets using a lexicographic order, each facet introduces at most one new
critical cell, while facet attachments leaving the homotopy
type unchanged do not contribute any new
critical cells.  The most natural of
poset edge-labellings and chain-labellings seem to yield lexicographic
discrete Morse functions in which one may easily give a systematic 
description of critical cells.  However,  
these are usually not minimal Morse functions.

The purpose of the present paper is to provide tools for 
``optimizing'' discrete Morse functions by cancelling pairs of critical 
cells, both in general and specifically in lexicographic discrete Morse
functions.  This paper develops and justifies machinery that has
recently been applied to examples in [12], [13], and as provided in 
Section 9.
We do not address the very interesting question of how to turn
an arbitrary discrete Morse function into an optimal one, but rather 
provide tools that seem to work well at turning very natural discrete
Morse functions on complexes that arise in practice,
especially on poset order complexes,
into ones with smaller Morse numbers (and often into optimal Morse 
functions).  These tools have led to connectivity 
lower bounds (in [13]) and to proofs that posets 
are Cohen-Macaulay, without needing to find a shelling (in [12]).  
Using our machinery typically seems to require somewhat lengthy proofs,
but this is because there are several 
things one needs to check, each of which is often straightforward.

We introduce and justify two new tools for cancelling critical
cells in discrete Morse functions, in Sections 3 and 7, respectively:
\begin{enumerate}
\item
A criterion for 
reversing several gradient paths,
each of which is individually reversible, to simultaneously cancel 
several pairs of critical cells.  
\item
A result showing that a pair of critical cells $\tau ,\sigma $
in a lexicographic discrete Morse function may be cancelled
whenever there is a certain type of gradient path from $\tau $ to $\sigma $,
by virtue of a correspondence between gradient paths and certain types of
reduced expressions for permutations
\end{enumerate}
Section 7 also introduces the notion
of a least-content-increasing labelling, a class of labellings for which
critical cell cancellation seems to be particularly manageable.
Section 4 very briefly shows how filtrations of
simplicial complexes may help in constructing discrete Morse
functions.  Section 5 reviews the notion of lexicographic
discrete Morse function, in preparation for later sections.  It is
quite helpful in cancelling critical cells 
to view these as coming from a filtration.
Section 6 gives a new characterization for 
lexicographically first reduced
expressions for permutations, as needed to prove (2) above.  

Finally, Sections 8 and 9 apply the above 
results to prove that two posets are homotopically
Cohen-Macaulay: the poset $PD(1^n,q)$
of collections of independent lines in a finite vector space,
and the poset $\Pi_{S_n}$ 
of partitions of 
$\{ 1,\dots ,n\} $ into cycles.  
The former is the face poset of a matroid complex, and so is 
well-known to be shellable (see [5]), but is included to illustrate on a 
familiar example a strategy for constructing and optimizing a discrete Morse 
function that has also worked on much more complex examples.
The latter poset was recently defined by
Remmel as a new partial order on elements of the symmetric group.  
See [12] and [13] for much more difficult applications of the machinery
provided in this paper.

\section{Background}

Recall that any permutation may be expressed as a product of adjacent
transpositions.  Denote by
$s_i$ the adjacent transposition swapping $i$ and
$i+1$.  An {\bf inversion} in a permutation $\pi $ is a pair $(i,j)$ with
$i<j, \pi (i)>\pi (j)$.
A {\bf reduced expression} for $\pi $ is a minimal expression for $\pi
$ as a product of adjacent transpositions.
Any two reduced expressions for the same permutation
are connected by a series of braid relations, i.e. relations 
$s_i\circ s_j = s_j\circ s_i$ for $|j-i|>1$ and $s_i\circ s_{i+1}
\circ s_i = s_{i+1}\circ s_i \circ s_{i+1}$.  See [Hum] or
[Ga] for more about reduced expressions for permutations.  
In our setting, permutations in $S_n$ act by permuting
positions of labels in a label sequence of length $n$, i.e. $s_i$ 
swaps the $i$-th and $(i+1)$-st labels.  

\subsection{Discrete Morse theory}
Forman introduced discrete Morse theory for CW complexes in
[10], then Chari provided a combinatorial reformulation for regular
cell complexes in [7].  Most discrete Morse functions in the literature
are in terms of this
combinatorial reformulation, in which one constructs an ``acyclic matching''
on the face poset of the complex.  

Let $\sigma^{(d)}$ denote a cell of dimension $d$.
Forman defines a function $f$ which assigns real 
values to the cells in a regular cell complex (as defined next)
to be a {\bf discrete Morse function} if 
for each cell $\sigma^{(d)}$, the
sets $ \{ \tau^{(d-1)} \subseteq \sigma^{(d)} | f(\tau^{(d-1)} ) 
\ge f(\sigma^{(d)} )\} $ and $ \{ \tau^{(d+1)} \supseteq \sigma^{(d)}
| f(\tau^{(d+1)} ) \le
f(\sigma^{(d)} )\} $ each have cardinality at most one.
Making these requirements for all cells implies that 
for each cell $\sigma $, at most one of the two cardinalities is nonzero.
When both cardinalities are 0, $\sigma $ is called a {\bf critical cell}.
For non-regular CW-complexes, Forman makes the additional
requirement $f(\sigma ) < f(\tau )$ 
for each $\sigma $ which is a non-regular face of $\tau $.  

Recall that
any CW complex $M$ is equipped with characteristic maps 
$h_{\tau }: B\rightarrow M$ sending a closed ball $B$ of dimension $p+1$ 
to $\overline{\tau^{(p+1)}}$.  
A face $\sigma^{(p)}$ is a {\bf regular face} of $\tau $ if the 
restriction of $h$ to $h^{-1}(\sigma )$ is a homeomorphism and the 
closure $\overline{h^{-1}(\sigma )}$ is a closed $p$-ball.  If all face
incidences in a CW complex are regular, then the complex is a {\bf regular
CW complex}.  All simplicial complexes and boolean cell 
complexes are regular (see [4], [21]).  
On the other hand, the minimal CW complexes for a 2-sphere
and the real projective plane have non-regular faces.

\begin{figure}[h]
\begin{picture}(250,100)(-60,15)
\includegraphics[width=0.4\textwidth, angle = 0]{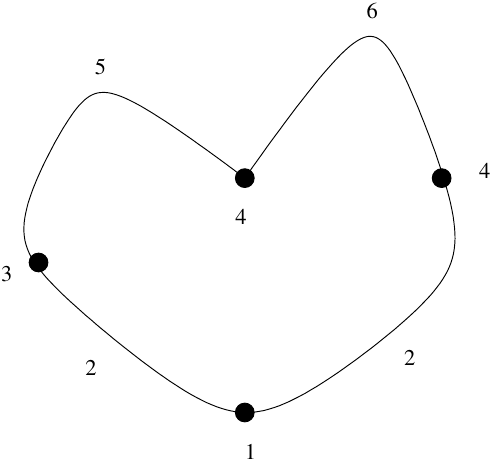}
\end{picture}
\caption{A discrete Morse function}
\label{height}
\end{figure}

Figure ~\ref{height} 
gives an example of a discrete Morse function on a 1-sphere resulting from
a height function.
Critical cells record changes in topological structure as a complex is 
built by sequentially inserting cells in the order specified by the Morse
function; non-critical cells may be eliminated by elementary
collapses without changing the homotopy type.
The non-critical cells come in pairs 
$\sigma^{(p)} \subseteq \tau^{(p+1)} $ which prevent each other from being 
critical by satisfying $f(\sigma )\ge f(\tau )$.  Elementary
collapses eliminating these pairs are possible because 
$\sigma $ is a free face of $\tau $ in a certain partial 
complex (see Section ~\ref{gradpath} for details).  

Any discrete Morse function on a regular cell complex gives rise to 
a matching on its face poset, by the aforementioned pairing on 
non-critical cells.  Recall that the {\bf face poset} 
of a regular cell complex is the partial
order on cells with
$\sigma < \tau $ for each $\sigma $ in the boundary of 
$\overline{\tau }$. 

\begin{defn}[Chari]
A matching on (the Hasse diagram of) the face poset of a regular cell complex
is {\bf acyclic} if the directed 
graph obtained by directing matching edges
upward and all other edges downward has no directed 
cycles.  Recall that the Hasse diagram of a poset is the graph
whose vertices are poset elements and whose edges are covering
relations, i.e. comparabilities $u<w$ such that there is no 
intermediate element $v$ satisfying $u<v<w$.
\end{defn}  
Observe that the face poset 
matching resulting from a discrete Morse function is always
acyclic, because the edges are oriented
in the direction in which $f$ decreases.  
Conversely,
many different (but in some sense equivalent) discrete Morse functions 
may be constructed from any face poset acyclic
matching.  

\begin{rk}
Given an acyclic matching on a facet poset with $n$ elements,
the corresponding discrete Morse functions are the points in a cone in
$\reals^n $ bounded by hyperplanes $x_i=x_j$ coming from pairs 
of comparable poset elements $v_i < v_j$.  The acyclic matching 
determines which side of each hyperplane contributes to the cone, and its
acyclicity ensures that the cone is non-empty.
\end{rk}

\begin{figure}[h]
\begin{picture}(250,60)(-72,10)
\includegraphics[width=0.3\textwidth, angle = 0]{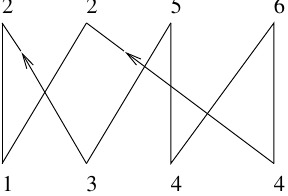}
\end{picture}
\caption{An acyclic matching}
\label{acyc_match}
\end{figure}

See Figure ~\ref{acyc_match} for the acyclic
matching corresponding to the discrete Morse function of Figure 
~\ref{height}.
 
Let $m_i$ be the number of critical cells of dimension $i$ in a discrete
Morse function on $\Delta $, and 
let $b_i$ denote the Betti number recording the rank of 
$H_i(\Delta )$.  
By convention, we include the empty set in our face
poset, so as to get a reduced version of discrete Morse theory, and we let
$\tilde{m}_i, \tilde{b}_i$ denote the reduced Morse numbers and reduced
Betti numbers, respectively.
Forman showed in [10] that $\Delta $ collapses onto a 
CW-complex $\Delta^M$ which has $m_i$ cells of dimension $i$ for each $i$,
such that $\Delta^M$ is homotopy equivalent to $\Delta $.  
The existence of such a complex $\Delta^M$ implies that
the following results from traditional Morse theory,
the first two of which are called the Morse inequalities, carry over to
discrete Morse theory.  Critical cells of dimension $i$ will play the 
role of critical points of index $i$.
\begin{enumerate}
\item
$\tilde{m}_j \ge \tilde{b}_j $ for $-1\le j\le dim(\Delta )$
\item
$\sum_{i=0}^{j+1}
(-1)^i \tilde{m}_{j-i} \ge \sum_{i=0}^{j+1} (-1)^i \tilde{b}_{j-i}$ 
for $0\le j\le dim(\Delta )$, with 
equality achieved when $j=dim(\Delta )$
\item
If $\tilde{m}_i=0$ for all $i\ge 0$, then $\Delta $ is collapsible.
\item
If $\tilde{m}_i=0$ for all $i\ne j$ for some 
fixed $j$, then $\Delta $ is 
homotopy equivalent to a wedge of $j$-spheres.
\end{enumerate}

\begin{rk}
Another consequence of the homotopy equivalence of $\Delta $ to such a 
complex $\Delta^M$ is that any
$\Delta $ with $\tilde{m}_i(\Delta )=0$ for all $i<j$ is $(j-1)$-connected.
\end{rk}

\begin{qn}
Is there a notion of
rank-selected lexicographic discrete Morse functions 
for graded posets?
\end{qn}

Figure ~\ref{height}
gives an example of a discrete Morse function with
$b_0=b_1 =1$ and $m_0=m_1=2$.  Letting $f(\emptyset )= 1.5$ turns this into
a reduced discrete Morse function with $\tilde{b}_0=0, \tilde{m}_0=
\tilde{b}_1=1$ and $\tilde{m}_1=2$.  The Morse 
numbers are larger than the Betti numbers because there is a critical cell
of dimension 0 that is labelled 4 which locally looks as 
though it is creating a new connected
component as the complex is built from bottom to top and there is a 
critical cell of dimension one that is labelled 5 which locally appears 
to be 
closing off a 1-cycle, but these two critical cells actually cancel each
other's effect.  

Recall that the {\bf order complex} of a finite poset $P$ with 
$\hat{0} $ and $\hat{1}$ is the simplicial complex, denoted 
$\Delta (P)$, whose $i$-faces
are chains $\hat{0} < v_0 < \cdots < v_i < \hat{1} $ of comparable poset 
elements.  Since $\tilde{\chi }(\Delta (P)) = \mu_P (\hat{0},\hat{1})$ (cf.
[19]), a discrete Morse function on $\Delta (P)$
gives a M\"obius function expression $\mu_P (\hat{0},\hat{1}) = 
\sum_{i=-1}^{\rm{dim}(\Delta(P))} 
(-1)^i \tilde{m}_i$, one of the original motivations
of [2].  A poset is {\bf homotopically Cohen-Macaulay} if each interval
is homotopy-equivalent to a wedge of spheres of top dimension.  This 
implies Cohen-Macaulayness of $\Delta (P)$ over any field, and for
$\rm{dim} (\Delta (P))>1$ that 
$\Delta (P)$ is simply-connected.  

\begin{rk}
\rm
{
A shelling (see [6] for background)
immediately implies the existence of a discrete Morse function
whose (reduced) critical cells are all top-dimensional (cf. [2], [7]).  
However, discrete
Morse functions may also give information about simplicial complexes and
CW-complexes that are far from shellable.
}
\end{rk}

Forman shows in [10] that whenever a discrete Morse function
has two critical cells $\sigma^{(p)} $ and $\tau^{(p+1)}$ such 
that there is a unique gradient path from $\tau $ to $\sigma $ (i.e. a path
upon which $f$ decreases at each step), then one
obtains a new acyclic matching in 
which $\sigma $ and $\tau $ are no longer 
critical by reversing this gradient path, analogously to in traditional
Morse theory.  
We call this process {\bf critical cell cancellation}.  
Uniqueness of the gradient path from $\tau $ to $\sigma $ implies 
that reversing it does not introduce any directed cycles.
We still get a matching because vertices along the 
path are matched with others along the path, and reversal
redistributes which pairs are matched so as
to incorporate the endpoints into the matching.  This reversal process 
for instance straightens the 1-sphere in Figure ~\ref{height}
into one in which two of the critical cells have been eliminated.  On the 
other hand, the minimal Morse numbers for the dunce cap are strictly larger
than its Betti numbers, reflecting the fact that it is contractible but not
collapsible.  

Section 7 provides machinery for checking uniqueness of a gradient path
in a lexicographic discrete Morse function.  In practice, we often
need to reverse several gradient paths simultaneously; 
section 3 introduces a notion of face poset for
$\Delta^M$, denoted $P^M$, and uses this
to provide a criterion for checking whether
several pairs of critical cells may be cancelled simultaneously.
Chari's construction does not apply directly 
because it is not clear a priori what are the face incidences in
$\Delta^M$, or 
how to define a face poset for a non-regular CW-complex:
for instance, one face may be incident to another in multiple ways (e.g. the
real projective plane, viewed as a CW-complex with one 0-cell, one 1-cell
and one 2-cell), or a cell may differ by more than one in dimension
from maximal cells in its boundary (e.g. a 2-sphere, realized as a
CW-complex with a 0-cell and a 2-cell).  

\section{Critical cell cancellation via 
multi-graph face posets}\label{gradpath}

Forman showed in [10] that a discrete
Morse function gives a way of collapsing a regular cell complex 
$\Delta $ onto a CW complex $\Delta^M $ of critical cells such that 
$\Delta \simeq \Delta^M$.  Theorem ~\ref{acyclic} 
follows from an analysis of the relationship between
face incidences in one such $\Delta^M$ and gradient paths in $\Delta $.  
First we recall from [17, Theorem 3.2]
one very explicit way of constructing such a complex $\Delta^M$, and
then make some observations about $\Delta^M$.

{\bf Constructing $\Delta^M$}: Begin with the empty complex and sequentially
attach cells in a way that illustrates how
$\Delta $ collapses onto $\Delta^M$, as follows.
At each step, attach a single cell not yet in the complex, all 
of whose faces have already been attached.  Choose any such 
cell $c_1$, and if it is critical then glue it to the 
complex.  Otherwise, $c_1$ is matched
with a cell $c_1'$, and by virtue of our insertion procedure,
$c_1$ will be in the boundary of $c_1'$.
If $c_1$ is the only cell in the boundary of $c_1'$ that
has not yet been inserted,
then we may attach the cells $c_1,c_1'$ to the complex both at once
without changing the homotopy type.
This is because $c_1$ is a free face
of $c_1'$, enabling an elementary collapse to eliminate $c_1, c_1'$ as
we construct $\Delta^M$.  
Kozlov used the acyclicity of the matching to show that there is an
attachment order for cells which allows each 
matched pair of non-critical cells to be inserted both at once with
one a free face of the other, and hence collapsed away.

\begin{rk} A critical cell $\sigma^{(p)} $ is incident to a 
critical cell $\tau^{(p+1)}$ in $\Delta^M $ if there
is a gradient path from $\tau $ to $\sigma $ in the discrete Morse 
function on $\Delta $.  Each such gradient path gives a distinct way
in which $\sigma $ is incident to $\tau $ in $\Delta^M$, 
implying that $\sigma $ is not a 
regular face of $\tau $ unless there is a 
unique gradient path from $\tau $ to $\sigma $.  
\end{rk}

An irregular face occurs in $\Delta^M$, for example, if two different 
$p$-faces of $\tau^{(p+1)}$ in $\Delta $ get identified
with the same face $\sigma $ through a series of collapses
leading from $\Delta $ to $\Delta^M$.  Alternatively, if there is 
a directed path in $F(\Delta )$ from a critical
cell $\tau^{(p+1)}$ to a critical cell $\sigma^{(q)}$ for 
$q < p$, and this path does not pass through any critical cells of 
intermediate dimension, then this also yields a face
$\sigma $ which is an irregular face of $\tau $ in $\Delta^M$. 
With these observations in mind, we define the {\bf multi-graph face poset},
denoted $P^M$, for the complex $\Delta^M$ of critical cells as follows:
\begin{enumerate}
\item
The vertices in $P^M$ are the cells in $\Delta^M$, or equivalently the 
critical cells in the discrete Morse function $M$ on $\Delta $.
\item
There is one edge between a pair of cells $\sigma^{(p)} ,\tau^{(p+1)}$ of 
consecutive dimension for each gradient path from $\tau $ to $\sigma $.
\end{enumerate}

\begin{figure}[h]
\begin{picture}(250,75)(0,0)
\includegraphics[width=0.7\textwidth, angle = 0]{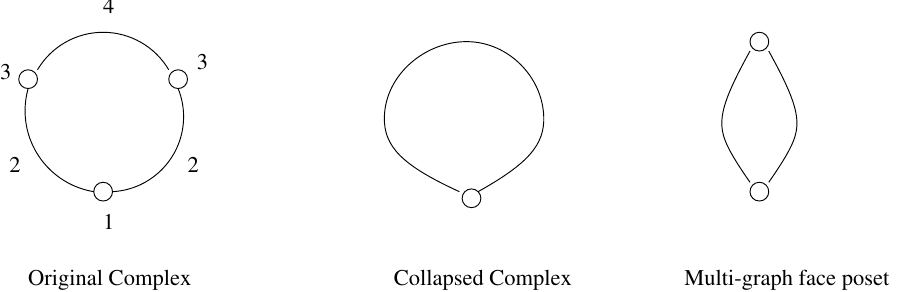}
\end{picture}
\caption{A ``face poset'' $P^M$ for $\Delta^M$}
\label{collapse}
\end{figure}

\begin{figure}[h]
\begin{picture}(125,50)(-35,10)
\includegraphics[width=0.15\textwidth, angle = 0]{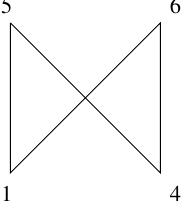}
\end{picture}
\caption{A multi-graph face poset}
\label{p_m}
\end{figure}

Figure ~\ref{collapse} 
gives an example of a discrete Morse function $M$ on a regular cell
complex $\Delta $, together with the non-regular CW complex $\Delta^M $ onto 
which $\Delta $ collapses and its multi-graph face poset $P^M$.
Recall that the discrete Morse function in Figure ~\ref{height} was
optimizable by reversing a single gradient path.  
Figure ~\ref{p_m} shows its multi-graph face poset, with critical cells
labelled by their Morse function values.  Note that any one of the 
four gradient paths between pairs of critical cells may be reversed to 
obtain an optimal Morse function.

For $f$ to be a discrete Morse function, 
Forman required $f(\sigma ) < f(\tau )$ whenever 
$\sigma $ is a non-regular face of $\tau $.
Notice that any discrete Morse function on 
$\Delta^M$ gives rise to an acyclic matching on $P^M$, by again
orienting edges in the direction in which $f$ decreases; 
we can never match faces $\sigma^{(p)},
\tau^{(p+1)}$ which are incident in multiple ways or which differ in 
dimension by more than one.  Next we verify
the converse, that any acyclic 
matching on $P^M$ implies the existence of a 
discrete Morse function on $\Delta^M$ whose critical cells are the 
unmatched elements of $P^M$.

\begin{prop}\label{rev_many}
Let $M$ be a discrete Morse function such that
there is a unique gradient path $\gamma_i $ from critical cell
$\tau_i $ of dimension $d_i$ to critical cell $\sigma_i $ of dimension
$d_i - 1$ for $1\le i \le r$. 
If there are no permutations $\pi\in S_r$ other than 
the identity such that there is a gradient path from 
$\tau_i $ to $\sigma_{\pi(i)}$ for $1\le i\le r$, then reversing the 
gradient paths $\gamma_i $ will not create any directed cycles.  
\end{prop}

\proof
Let us initially assume that the gradient
paths are non-overlapping.  When $r=2$, we can reverse 
both gradient paths $\gamma_1,\gamma_2$ without creating any directed 
cycles unless there is a cycle involving the reversals of both 
$\gamma_1 $ and $\gamma_2$.
The existence of a cycle involving $\gamma_1^{rev}$ and $\gamma_2^{rev}$
would imply that
$\gamma_i $ has elements $v_{i,1},v_{i,2}$ for $i=1,2$ such that the original 
digraph has paths from from $v_{1,1}$ to $v_{2,2}$, from $v_{2,1}$ 
to $v_{1,2}$, and from $v_{i,1}$ to $v_{i,2}$ for $i=1,2$.  
Thus, we obtain directed paths
$\tau_1\rightarrow v_{1,1}\rightarrow v_{2,2}\rightarrow\sigma_2$
and $\tau_2\rightarrow v_{2,1}\rightarrow v_{1,2}\rightarrow\sigma_1$
in the original directed graph. 
It is straightforward to extend this to larger $r$.

Now suppose some edge is shared by two different gradient paths 
$\gamma_i $ and $\gamma_j$.  Then the first and
last edges shared by $\gamma_i ,\gamma_j$ must be oriented upward, 
so that their endpoints are each only 
in a single matching edge.  We get gradient paths
from $\tau_i $ to $\sigma_j $ and $\tau_j$ to $\sigma_i $ by switching
gradient paths at the shared edge, implying the 
existence of a transposition $\pi $ of the type that is forbidden.
\EOP

\begin{rk}
{\rm
An easy way to ensure that a 
collection of gradient paths $\gamma_1,\dots ,\gamma_r $ which are 
individually reversible will also be simultaneously reversible
is for the discrete 
Morse function $f$ to satisfy $f(\tau_i ) < f(\sigma_j)$ for each pair
$i<j$.  This was the approach taken in [2].
}
\end{rk}

Subsequent applications (in [12], [13] and later sections of
this paper) seem to require a more careful analysis, as provided
in Theorem ~\ref{acyclic}.

\begin{thm}\label{acyclic}
Any acyclic matching on $P^M$ gives rise to a discrete Morse function for
$\Delta $ in which the critical cells of $\Delta $ are the unmatched elements
of $P^M$.  
\end{thm}

\proof
Let $M$ be the acyclic matching on $F(\Delta )$ that gave rise to
$\Delta^M $ and  $P^M$.  We must show that
simultaneously reversing all the gradient
paths in $F(\Delta )$ which correspond to edges in an acyclic matching for
$P^M$ gives a new acyclic matching on $F(\Delta )$.
There are three things to check, the third of which was already done
in Proposition ~\ref{rev_many}.
\begin{enumerate}
\item
that reversing this set of gradient paths is a well-defined operation,
i.e. there
are no edges in the Hasse diagram for $F(\Delta )$ that must be reversed
in one gradient path and cannot be reversed in another.
\item
that reversing these gradient paths still gives a matching.   
\item
that reversing all of these gradient paths does not create any directed 
cycles.
\end{enumerate}

To prove 1, note that any shared gradient path segments must be of the form 
$UD\dots DU$ so as to come from a matching, where 
$U$ and $D$ denote upward and downward-oriented
edges respectively.  If one gradient path $\gamma_1$
needs this segment reversed, while another gradient path $\gamma_2$
containing the segment is not reversed, we
reverse the shared segment without reversing all of $\gamma_2$, to
nonetheless get a matching.  Acyclicity is preserved  
by similar reasoning to that used in Proposition ~\ref{rev_many}.

The second requirement follows from the fact that vertices along a
gradient path are matched with others
along that gradient path.
The only potential difficulty would be 
if two gradient paths shared a segment $UD\dots DU$ and 
we reversed both
paths, but then the matching for $P^M$ would have a directed cycle.
\EOP

Section ~\ref{lex} will provide a way of determining when
a gradient path between a single pair of critical cells in a 
lexicographic discrete Morse function
is reversible.   This result and Theorem ~\ref{acyclic}
appear to be most useful when applied jointly.

\section{Discrete Morse functions on filtrations}\label{filtration}

This section shows how a filtration may 
simplify the task of constructing an acyclic matching by splitting it into
smaller, typically much more manageable pieces.  We independently 
discovered Lemma ~\ref{filter}, but it also appears as the ``Cluster
Lemma'' in [16], and the idea has been used widely 
(e.g. see [1], [14], [17], and [20]).  
The argument below also applies
to non-regular CW complexes obtained from regular ones by a series of 
collapses, using the generalized notion of acyclic matching from Section 
~\ref{gradpath}.  
\begin{lem}\label{filter}
Let $\Delta $ be a regular $CW$ complex which decomposes into collections 
$\Delta_{\sigma }$ of cells, indexed by the elements $\sigma $ in a partial 
order $P$ which has a 
unique minimal element $\hat{0}=\Delta_0$.  Furthermore, assume
that this decomposition is as follows:
\begin{enumerate}
\item
Each cell belongs to exactly one $\Delta_{\sigma }$ .
\item
For each $\sigma\in P$, $\cup_{\tau \le \sigma} \Delta_{\tau}$ 
is a subcomplex of $\Delta $.
\end{enumerate}
For each $\sigma\in P$, let 
$M_{\sigma }$ be an acyclic matching on the subposet of $F(\Delta )$
consisting of the cells in $\Delta_{\sigma }$.
Then $\cup_{\sigma\in P} M_{\sigma }$ is an 
acyclic matching on $F (\Delta )$.
\end{lem}

\proof
Let $D$ be the directed graph obtained by orienting matching edges 
upward and all other edges in $F(\Delta )$ downward.  By design,
$D$ does not include any upward-oriented edges between cells in different
components.  Suppose there is a downward-oriented edge from a cell in 
$\Delta_{\tau }$ to a cell in $\Delta_{\sigma }$ for $\sigma \ne \tau $.  
This
implies $\sigma \le \tau $, because $\cup_{\rho\le \tau} \Delta_{\rho }$ 
is a 
subcomplex of $\Delta $.  If $\sigma
\ne \tau $, then $\sigma < \tau $, which means it will be impossible for 
the directed path ever to return to $\Delta_{\tau }$.
\EOP

\medskip
The lexicographic discrete Morse functions of [2] to be discussed 
shortly may be viewed as coming from a filtration
$\Delta_1 \subseteq \cdots \subseteq \Delta_k$ with
$\Delta_j \setminus \Delta_{j-1} = F_j \setminus (\cup_{i<j} F_i)$, where
$F_1,\dots ,F_k$ is a lexicographic order on the facets of an order
complex.  

\section{Review of lexicographic discrete Morse functions}\label{lexrev}

This section briefly reviews lexicographic discrete Morse functions 
for poset order complexes from [2], 
in preparation for later sections.  We also briefly indicate how 
lexicographic discrete Morse functions yield an easy new proof 
that every interval in the weak order for $S_n$ is a 
ball or a sphere of specified dimension.

Any lexicographic order on the saturated chains of any finite poset $P$ 
with minimal and maximal elements $\hat{0}$ and $\hat{1}$ gives rise to a 
discrete Morse function on its order complex 
with a relatively small number of critical cells as follows.
Let $\lambda $ be a labelling of Hasse diagram edges (or more generally,
a chain-labelling) with integers such
that $\lambda (u,v) \ne \lambda (u,w)$ for $v\ne w$.  We obtain a total
order $F_1,\dots ,F_r$ on facets of $\Delta (P)$ by lexicographically
ordering the label sequences on saturated chains.  

\begin{rk}
Often we will
refer to ranks of elements in a maximal chain.  We do not assume
posets are graded, but rather allow the rank of an element
to depend on the choice of maximal chain within which it is considered.
\end{rk}

A key observation for
the [2] construction was
that each of the maximal faces in $F_j \cap (\cup_{i<j} F_i)$ has
rank set of the form, $1,\dots ,i,j,\dots ,n$, i.e. it omits a 
single interval $i+1,\dots ,j-1$ of consecutive ranks, by virtue of
our use of a lexicographic order.  Call the rank interval $[i+1,j-1]$
a {\bf minimal skipped interval} of $F_j$ with {\bf support} 
$i+1,\dots ,j-1$ and {\bf height} $j-i-1$.  Call the collection of 
minimal skipped intervals for $F_j$ the {\bf interval system} of $F_j$.
  
Notice that the faces 
in $F_j \setminus (\cup_{i<j} F_i )$, i.e. belonging to $F_j$ but not to
any earlier facets, are the faces of $F_j$
that include at least one rank from each of the minimal skipped intervals
of $F_j$.  For each $j$, [2] constructs an acyclic matching on the
set of faces in $F_j \setminus \cup_{i<j} F_i$.  It is immediate from
Proposition ~\ref{filter} (and was verified by other means in [2])  
that the union of these matchings is acyclic on $\Delta (P)$. 

Each $F_j \setminus (\cup_{i<j} F_i)$ has
a single critical cell if the homotopy type of $\Delta (P)$
changes with the attachment of $F_j$, and otherwise  
$F_j \setminus (\cup_{i<j} F_i)$ contains no
critical cells.  (See [3] for a related notion of 
shellability, called a weak shelling.)  We say 
that $F_j$ {\bf contributes} a critical cell when
$F_j\setminus (\cup_{i<j} F_i)$ has a critical cell.  
$F_j$ contributes a critical cell if and only if the interval system
of $F_j$ covers all ranks in $F_j$ after the truncation procedure 
described shortly; the dimension of the critical cell
is then one less than the number of intervals in the truncated 
interval system.

For convenience,
order the minimal skipped intervals $I_1,\dots ,I_m$ of $F_j$ so
that their lowest rank elements sequentially increase in rank.  

\medskip
{\bf Description of critical cells in a lexicographic discrete Morse function:} 
\begin{itemize}
\item
if the minimal skipped intervals of $F_j$ (namely the $I$-intervals) do 
not collectively have support covering
all the ranks in $F_j$, then $F_j$ does not contribute a critical cell
\item
if the $I$-intervals of $F_j$ 
have disjoint support covering all ranks in $F_j$, then the critical
cell in $F_j \setminus (\cup_{i<j} F_i)$ consists of
the lowest rank from each of the minimal skipped intervals
\item
if there is overlap in the minimal skipped intervals of $F_j$, but they
cover all ranks in $F_j$, then iterate the following procedure to 
obtain a potential critical cell:
\begin{enumerate}
\item
include the lowest rank from $I_1$ in the critical cell
\item
truncate all the remaining minimal 
skipped intervals by chopping off any ranks that they share with $I_1$
\item
discard $I_1$ and any skipped intervals that now strictly contain other
intervals in our collection
\item
re-index the remaining truncated minimal skipped 
intervals to begin with a newly chosen $I_1$
\item
repeat until there are no more minimal skipped intervals
\end{enumerate}
\end{itemize}

\begin{defn}
The truncated, minimal intervals obtained by the above 
procedure are called the $J$-intervals of $F_j$, and are
nonoverlapping.  If the $J$-intervals do not cover all ranks, then
$F_j$ does not contribute a critical cell.
\end{defn}

\begin{rk}
Since each maximal chain contributes at most one critical cell, we 
refer interchangeably 
to critical cells and to the maximal chains contributing them.
For instance, the label sequence of a critical cell will refer to 
the label sequence for the maximal chain which contributes the critical 
cell.
\end{rk}

To construct a lexicographic discrete Morse function and compute its
Morse numbers requires a labelling in which one may understand its 
minimal skipped intervals, so as to determine 
which saturated chains contribute critical cells and what are their
dimensions.  
To cancel critical cells by gradient path reversal, one
also needs to understand the matching construction of [2] 
enough to recognize gradient paths between critical cells and check their
uniqueness.  

\medskip
{\bf Description of acyclic matching on $F_j \setminus \cup_{i<j} F_i$:}
\begin{itemize}
\item
If the $I$-intervals leave some rank uncovered, then there is a
cone point in $F_j \cap (\cup_{i<j} F_i)$, so match by 
including/excluding the lowest rank such cone point.  
\item
Otherwise, match any cell that differs from the potential critical
cell on some $J$ interval based on 
the lowest rank $J$-interval where it 
differs from $F_j$'s potential critical cell.  Specifically, 
match by including/excluding the element of lowest rank in this 
interval.
\item
If the $I$-intervals cover all ranks but the $J$-intervals leave 
some rank uncovered, match all remaining cells by including/excluding 
the element at the lowest such uncovered rank.
\end{itemize}

\begin{rk}
\rm{
If all the $I$-intervals in a lexicographic discrete Morse function on a 
poset $P$ have 
height at most $d-1$, then 
each interval $(x,y)$ in $P$ is at least
$(-1 + \frac{\rm{rk}(y) 
- \rm{rk}(x) -1}{d-1} )$-connected.  
}
\end{rk}

\begin{examp}\label{lex-examp}
\rm{
Consider the ring $k[ab,a^2,c,d,e,b^2] \cong 
k[x_1,\dots ,x_6]/(x_2x_6 - x_1^2)$ 
and the partial order on (equivalence classes of) monomials by 
divisibility.  Label $m_1\prec m_2$ by the 
quotient $m_2/m_1$.  Order 
the labels $x_i < x_j$ for each $i<j$.  
Consider the interval $(1,x_2x_3x_4x_5x_6)$, and its 
saturated chain $F_j = 1\prec x_2 \prec x_2x_3\prec x_2x_3x_6
\prec x_2x_3x_5x_6 \prec x_2x_3x_4x_5x_6$, depicted in
Figure ~\ref{critcell}.  
This chain is labelled $x_2x_3x_6x_5x_4$, or more compactly
$23654$ (recording indices).

\begin{figure}[h]
\begin{picture}(250,150)(-60,10)
\includegraphics[width=0.4\textwidth, angle = 0]{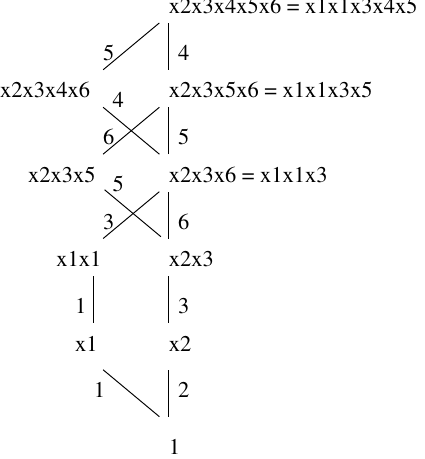}
\end{picture}
\caption{A saturated chain and its minimal skipped intervals}
\label{critcell}
\end{figure}

$F_j$ has minimal skipped intervals $I_1, I_2, I_3$, with
support $\{ 1,2\} $,
$\{ 3\} , $ and $\{ 4\} $, respectively.  $I_1$ results from
the syzygy $(x_2x_6 - x_1^2)x_3 =0$, causing a lexicographically
earlier facet $F_i$ labelled $11354$.
The intervals $I_2$ and $I_3$ come from descents.  The intervals
$I_1, I_2, I_3$ are non-overlapping and cover $F_j$, so $F_j$ contributes
a critical cell $Crit(F_j)$ with rank set $\{ 1, 3, 4\}$.
$I_1, I_2$, and $I_3 $.  Thus, $Crit(F_j) = x_2 < x_2x_3x_6 
< x_2x_3x_6x_5$.
}
\end{examp}

\subsection{Example: the weak order for the symmetric group}

In [9], Edelman and Walker study the poset of regions in a hyperplane
arrangement with a chosen base region, partially ordered by inclusion
of the set of hyperplanes separating a region from the base region $B$.
They show that each interval is a ball or a sphere of a 
proscribed dimension.
The weak order for any Coxeter group may be viewed as
the poset of regions for a Coxeter arrangement.  

\begin{rk}
The weak order for the symmetric group has a lexicographic discrete Morse
function yielding an elementary proof of the result from [9] in this case.
\end{rk}

The idea is to label each covering 
relation by the simple reflection being applied.  
The minimal skipped intervals in the resulting lexicographic
discrete Morse function are 
labelled by reduced expressions of the following forms:
\begin{enumerate}
\item
$s_j \circ s_i $ for $j>i+1$, since $s_i \circ s_j$ is lexicographically 
earlier
\item
$s_{j+1}\circ s_j \circ s_{j-1} \circ \cdots \circ s_{i+1} \circ s_i \circ
s_{j+1}$ for some $i\le j$, since $s_j \circ s_{j+1}\circ s_j \circ s_{j-1}
\circ \cdots \circ s_{i+1} \circ s_i$ is lexicographically earlier
\end{enumerate}

Within any interval it is straightforward
to cancel all but at most critical cell, using Theorems
~\ref{acyclic} and ~\ref{red_exp}, and to see that the dimension 
of any surviving critical cell agrees with the result of [9].  

\begin{qn}
May lexicographic discrete Morse functions be used to deduce
more general results about posets of regions or be applied to the (open)
question of whether intervals in higher Bruhat
orders are homotopy equivalent to balls or spheres?
\end{qn}

\section{Lexicographically first reduced expressions for permutations}
\label{lex-first-sec}

This section gives a new characterization for lexicographically first
reduced expressions for permutations, 
in the sense of [8].  This provides a very natural explanation for
[8, Theorem 2.5], which characterized 
the ``type'' of any lexicographically
first reduced expression, i.e. the possible vectors $(m_1,\dots ,m_{n-1})$
where $m_i$ counts the number of appearances of the simple reflection 
$s_i$ in a lexicographically first reduced expression.  
Our main interest in Theorem ~\ref{lex-last} below is that 
it will help us in Section 7 to show that certain gradient paths in
lexicographic discrete Morse functions are unique and hence reversible.

Recall that a reduced expression for a permutation $\pi $ is the 
{\bf lexicographically first reduced expression } for $\pi $ if its reduced
word precedes all other reduced words for reduced expressions for $\pi $ in
lexicographic order, i.e. when we say $w_1 < w_2$ if $w_1$ has a smaller 
letter at the first place that the two words differ.

\begin{thm}\label{lex-last}
Every permutation has a unique reduced expression 
$s_{i_1}\circ \cdots \circ s_{i_k}$ with 
all of the following properties:
\begin{enumerate}
\item
$s_i$ is never immediately followed by $s_j$ for $j>i+1$
\item
There are no reduced expressions in the commutation class for 
$s_{i_1}\circ \cdots \circ s_{i_k}$ which have 
$s_{i+1}\circ s_i\circ s_{i+1}$
appearing consecutively.
\end{enumerate}
\end{thm}

The proof is somewhat similar to a proof in [11] for the fact 
that any two reduced expressions for a fixed permutation are 
connected by a series of braid relations.  

\proof
Existence follows from the fact  
that the reversal of the lexicographically 
first reduced expression for the inverse permutation must take this form.
To prove uniqueness, first check that any allowable reduced 
expression $\gamma $
has exactly one copy of $s_u$ for $u$ the maximal index appearing in 
$\gamma $.
Otherwise there would be a 
copy of $s_{u-1}$ between each pair of copies of $s_u$;
avoiding $s_u\circ s_{u-1}\circ s_u$ in the 
commutation class would necessitate another $s_{u-1}$
between the same two copies of $s_u$, forcing an intermediate
$s_{u-2}$, and this continues indefinitely, contradicting
the finiteness of $\gamma $.

Similarly for each $t\le u$,
each pair of $s_{t-1}$'s in $\gamma $ is separated by 
an $s_t$.
Moreover, each $s_{t-1}$ which is followed by an $s_t$ must 
immediately precede an $s_t$.
If $\gamma $ includes
$s_t$ and of $s_{t-d}$ for $d>1$, but 
no copies of $s_{t-d'}$ for $1\le d' < d$,
then there is at most one copy of $s_{t-d}$ and it must appear after the
last copy of $s_t$.  
Thus, the only flexibility is in which
allowable adjacent transpositions do appear, but varying this
will change the permutation itself.  Thus, there is at most one allowable
reduced expression for any fixed permutation.
\EOP

\begin{cor}
A reduced expression $\omega = s_{i_1}\circ \cdots \circ s_{i_k}$
is the lexicographically first reduced expression for
a permutation $\pi $ if and only if all of the following hold: 
\begin{enumerate}
\item
$\omega $ never has $s_j$ immediately followed by $s_i$ for $i< j-1$
\item
$\omega $ never has consecutive simple reflections $s_{i+1}\circ s_i\circ 
s_{i+1}$
\item
$\omega $ is not in the same commutation class as any reduced expressions
containing consecutive simple reflections $s_{i+1}\circ s_i\circ s_{i+1}$ 
\end{enumerate}
\end{cor}
\proof
There is a unique reduced expression for $\pi $ satisfying these conditions,
by a nearly identical proof to the one used for Theorem ~\ref{lex-last}.  
These conditions are necessary, because 
otherwise a lexicographically earlier reduced expression
could be obtained by applying a Coxeter relation.  
\EOP

In the language of [8], the third condition says that $\omega $ is not
$C_1$-equivalent to any reduced expressions containing consecutive
simple reflections $s_{i+1}\circ s_i\circ s_{i+1}$.  

\begin{rk}
Theorem 2.5 of [8]
showed that a vector $(m_1,\dots ,m_{n-1})$ is the type vector of a 
lexicographically first reduced expression if and only
if $0\le m_{n-1}\le 1$ and $0\le m_i \le m_{i+1} + 1$ for each 
$1\le i\le n-2$.
This also follows immediately from the arguments used in the proof of 
Theorem ~\ref{lex-last} above.  
\end{rk}

\section{Gradient paths in lexicographic discrete Morse functions}\label{lex}

This section uses
properties of reduced expressions for permutations 
to give conditions under which a gradient path between a pair of
critical cells in a lexicographic discrete Morse function is unique,
and hence may be reversed to cancel the pair of 
critical cells.  While the results in this section are quite
technical,
they have recently been useful in applications (see [12] and [13],
as well as later sections of this paper).  We also believe
the proofs in this section may be generalizable, for instance to deal 
with nonsaturated chain segments, at least in special cases related to 
least-content-labellings, as defined shortly.

\begin{defn}
A poset edge labelling is {\bf least-increasing} if every interval has a
(weakly) increasing chain as its lexicographically smallest saturated chain.
It is {\bf least-content-increasing} if in addition the 
lexicographically smallest chain is at least as small as the 
increasing rearrangement of every other label sequence on the interval.
Within a least-increasing labelling, any 
increasing chain that is not lexicographically smallest on its 
interval is called a {\bf delinquent chain}.
\end{defn}

Later we will show that least-content-increasing labellings 
satisfy the fairly technical
conditions of Theorem ~\ref{red_exp} below.
Notice that the least-increasing condition is less restrictive than an 
EL-labelling, since intervals may
have several increasing chains.
See Sections ~\ref{pd1nq} and ~\ref{pisn} as
well as [13] for examples of least-content-increasing labellings.
Least-content-increasing labellings are particularly well-suited to 
cancelling pairs of critical cells whose label sequences have equal
content, i.e. equal unordered multiset of labels,
using the following observation:

\begin{rk}
If $\tau ,\sigma $ are critical cells whose label sequences have equal
content in a lexicographic discrete Morse function coming from
a least-content-increasing labelling, then any gradient path from 
$\tau $ to $\sigma $ must preserve content at each step and also may
never introduce inversions.  Thus, gradient path steps must sort
labels whenever the label sequence changes at all.
\end{rk}

Any gradient path from a critical cell $\tau^{(p+1)}$ to a critical
cell $\sigma^{(p)}$ must alternate between ranks $p+1$ and $p$, because it
cannot contain two consecutive upward steps or end with an upward step.
Denote the downward step in a gradient path which deletes the $i$-th element
from a chain by $d_i$.
Denote by $u_r$ the unique upward step from a non-critical cell to its 
matching partner, if the newly inserted chain element is the $r$-th
element in the chain.  We will focus on gradient paths 
resulting from labellings such that each
$d_i$ deletes from a chain a descent
and then $u_i$ replaces this by a lexicographically 
earlier ascent in which the two labels have been swapped.  Thus, we may
view the pair of steps $d_i \circ u_i$ as an adjacent transposition
$s_i = (i, i+1)$ swapping the labels in positions $i,i+1$ on the 
lexicographically earliest saturated chain containing the order complex 
face.

\begin{defn}
An adjacent transposition $s_i$ {\bf acts effectively} on a chain $C$ 
that includes 
elements $v_{i-1},v_i,v_{i+1}$ of ranks $i-1,i,i+1$ and has a descent
at rank $i$
if $s_i$ sends $C$ to a chain in which $v_i$ is replaced by 
$v_i'$ such that $\lambda (v_{i-1},v_i) = \lambda(v_i',v_{i+1})$ and 
$\lambda (v_{i-1},v_i') = \lambda (v_i,v_{i+1})$.
A reduced expression $\omega $ {\bf acts effectively} on $C$
if each of its adjacent transpositions acts effectively in turn.  
\end{defn}

Let $\lambda (u,v)$ denote the sequence of edge labels on the 
lexicographically
earliest saturated chain from $u$ to $v$.  Then $s_i$ {\bf acts 
effectively} on a (not necessarily saturated) chain $C = \hat{0} < 
v_1 < \cdots < v_{i-1}
< v_i < v_{i+1} < \cdots < v_{r-1} < \hat{1}$ 
if $s_i(C)$ replaces $v_i$ by some $v_i'$ 
such that $rk(v_i')=rk(v_{i-1})+1$ and the concatenated label sequence
$\lambda(v_{i-1},v_i')\circ \lambda(v_i',v_{i+1})$ is obtained from
$\lambda(v_{i-1},v_i)\circ \lambda(v_i,v_{i+1})$ by sorting 
labels into increasing order.

\begin{defn}
A poset edge-labelling (or chain-labelling) 
has the {\bf ordered zigzag property} if 
$\lambda (u_i,v_i)<\lambda (u_i,v_{i+1})$ for $1\le i\le k$ implies that
$\lambda (u',v_1)<\lambda (u',v_{k+1})$ for every $u'$ which satisfies
both $u'<v_1$ and $u'<v_{k+1}$.
\end{defn}

\medskip
Any product of chains has the ordered zigzag property when each 
covering relations is labelled by the coordinate being increased.
To apply Theorem ~\ref{red_exp}, it will suffice to have the ordered zigzag
property on a subposet which includes the two critical cells $\sigma ,\tau $ 
to be cancelled as well as all cells with intermediate Morse function 
values.  See [13] for such an application.

In the next theorem, we will use reduced expressions for permutations 
to show that certain gradient paths are reversible.
By convention, we apply adjacent transpositions from left to right.  
Let $\tau^{(p+1)}$ be a critical cell which includes
ranks $i_0,i,i+1,\dots ,j-1,j$.  Denote by $u_{\tau }, v_{\tau }$ the 
elements at ranks $i_0, j$, respectively in $\tau $, and say that these
are the $k$-th and $l$-th elements of $\tau $, respectively.  Let 
$\sigma^{(p)}$ be a critical cell that agrees with $\tau $ except between
ranks $i_0$ and $j$.  Theorem ~\ref{red_exp} will
show that when a certain type of gradient path from $\tau $ to $\sigma $
exists, then it is the only gradient path from $\tau $ to $\sigma $, 
and hence may be reversed to cancel $\tau $ and $\sigma $.

\begin{thm}\label{red_exp}
Let $M$ be a lexicographic discrete Morse function 
satisfying the ordered zigzag property, with critical cells $\tau ,\sigma $
as above.  Suppose for each $k< r < l$, 
$s_r$ acts effectively on any chain $C$ that (1) coincides
with $\tau $ except strictly between ranks $i_0$ and $j$,
(2) has a descent at rank $r$, and (3) includes ranks $r\pm 1$.
Let $\omega $ be a reduced expression 
that acts effectively on $\tau^{(p+1)}$ in such a way that $\omega $
followed by some $d_r$ for $k<r<l$ yields a critical cell
$\sigma^{(p)}$ that agrees with $\tau $ except between ranks $i_0$ and $j$. 
Then the gradient path from 
$\tau $ to $\sigma $ is unique unless it ends with steps 
$s_{r+1}\circ s_{r} \circ d_{r+1}$ 
or differs by commutation relations from
one ending this way, in which case there are 
exactly two gradient paths from $\tau $ to $\sigma $.
\end{thm}

\proof
Let $\omega^* = \omega \circ s_i$, where $i$ is the
rank of the final downward step $d_i$ in the given
gradient path from $\tau $ to $\sigma $.
We break the proof into the following three parts:
\begin{enumerate}
\item
If $\omega^* $ involves only
$s_{k+1},\dots ,s_{l-1}$,
then $\omega^* $ 
does not contain any consecutive adjacent transpositions
$s_t\circ s_u$ for $u>t+1$.  It also avoids consecutive adjacent 
transpositions
$s_{t+1}\circ s_t\circ s_{t+1}$, except when the second
$s_{t+1}$ occurs as the final downward
step.  Furthermore, $\omega^* $ is not equivalent up to commutation 
relations to any reduced 
expression involving $s_{t+1}\circ s_t\circ s_{t+1}$ except with the second
$s_{t+1}$ as the final downward step.

\item
Let $\pi $ be the permutation which has $\omega^*$ as one of its reduced
expressions.  Every 
gradient path $\gamma $ from $\tau $ to $\sigma $ corresponds to a 
reduced expression for $\pi $. 
That is, $\gamma $ has the form 
$d_{i_1}\circ u_{i_1}\circ \cdots \circ d_{i_q}$
for some indices $i_1,\dots ,i_q$ such that 
$s_{i_1}\circ \cdots \circ s_{i_q}$ is a 
reduced expression for $\pi$. 
\item
Every permutation has a unique reduced expression that avoids consecutive
transpositions $s_i\circ s_j$ for $j>i+1$ and whose commutation class 
avoids $s_{i+1}\circ s_i \circ s_{i+1}$.  In addition, every permutation
has at most one reduced expression that avoids consecutive adjacent
transpositions $s_i\circ s_j$ for $j>i+1$ and whose commutation class
avoids $s_{i+1}\circ s_i \circ s_{i+1}$ except as the final three steps.
Thus, the gradient path $\omega $
from $\tau $ to $\sigma $ is unique, except in the designated case 
where there are two gradient paths.
\end{enumerate}
We begin with (1).
The series of steps $d_t\circ u_t \circ d_u$ for $u>t+1$
yields a chain lacking rank $u$, and having an ascent at rank
$t$.  Thus, it is the top of an ``up'' edge in 
$F(\Delta (P))$; it is matched with
the lower-dimensional face in which the cone point coming from the ascent
at rank $t$ is deleted.  Thus, the gradient path cannot continue.
Similarly, the
steps $d_{t+1}\circ u_{t+1}\circ d_t\circ u_t\circ d_{t+1}$ yield a chain 
with rank $t+1$ omitted, having an ascent at rank $t$,
precluding continuation of the gradient path.  
A reduced expression in the same commutation class as one
containing $s_{t+1}\circ s_t\circ s_{t+1}$ also gets stuck after the 
second $d_{t+1}$, by the same reasoning, unless this is the final step.  

When $\omega^* $ is equivalent up to 
commutation to a reduced expression which 
ends with steps $s_{t+1}\circ s_t\circ s_{t+1}$, we obtain the second
gradient path as follows.  Apply the braid relation to obtain $s_t
\circ s_{t+1} \circ s_t$ then consider the unique reduced expression
obtained from this by commutation relations, which 
avoids ever having $s_i$ followed 
by $s_j$ for $j>i+1$.  This gradient path
will end with $d_t \circ u_t$.  
This does not give the desired critical cell because the wrong ranks are
covered on $J$-intervals, but we may apply $d_{t+1}$ to obtain the 
desired critical cell.

Now we turn to (2). 
Notice that $\tau $ and $\sigma $ must agree up through rank $i_0$, 
since $\omega $ leaves these ranks unchanged.  Thus,
every gradient path from $\tau $ to $\sigma $ must also leave 
ranks $1,\dots ,i_0$ fixed in order for the discrete Morse function $f$ 
to be non-increasing along the gradient path from $\tau $ to $\sigma
$.  The ordered zigzag 
property ensures that every gradient path from 
$\tau $ to $\sigma $ also must leave ranks $j$ and higher untouched,
so only uses $d_r,u_r$ for $k<r<l$.  Each such $d_r$ must be followed 
by the upward step $u_r$, because applying $d_r$ to a chain $C$ ensures
that $d_r(C)$ has an ascent (and hence a cone point in $d_r(C) \setminus
(\cup_{F_i <_{lex} d_r(C)} F_i )$) on the newly uncovered interval.
There must be no
lower cone points, or we would be stuck at the top of an up edge.  
Each pair of steps $d_r \circ u_r$ for $r>j+1$
swaps two consecutive chain labels if the labels
are out of order prior to $d_r$, so that 
a series of effective steps $d_{i_1}\circ u_{i_1}\circ
\cdots \circ d_{i_k}\circ u_{i_1}$ eliminates exactly
the same inversions that $s_{i_1}\circ \cdots\circ s_{i_q}$ would.
Thus, every gradient path corresponds to an expression
for $\pi $ as a product of adjacent transpositions.

What remains to show is that non-reduced expressions
do not give rise to gradient paths.  Such an  
expression will eventually apply some adjacent transposition
$s_r$ at a rank $r$ where the labels $\lambda_1, \lambda_2$ are
increasing.  The step $d_r$ would not go in the direction the Morse 
function decreases, because the edge would be a matching
edge if $\lambda_1 ,\lambda_2$ is the lowest increasing pair of 
consecutive labels in the chain obtained by deleting rank $r$; otherwise,
the subsequent upward step would not be a matching step,
because the chain would instead match by inserting or deleting a lower 
rank cone point.  Either way, the gradient path cannot be completed.

The first statement in Part (3) was proven as Theorem ~\ref{lex-last}.
The second statement is proven using similar ideas.  
The only variation needed is 
that once we may have two copies of $s_i$ after the final $s_{i+1}$, but
then we must have one copy of $s_{i-1}$ between the copies of $s_i$.  
The rest of the analysis proceeds just as in Theorem ~\ref{lex-last}.
\EOP

\begin{rk}\label{dual_red_exp}
Theorem ~\ref{red_exp} also applies to critical cells $\tau ,\sigma $
where $\tau $ includes ranks $i,i+1,\dots ,j-1,j,k$ for $k>j+1$,
$\sigma $ includes ranks $i,i+1,\dots ,j-1,k$, and $\tau, \sigma $
agree outside of the interval from rank $i$ to $k$.
\end{rk}

\begin{cor}
If the conditions of the previous theorem are met, and in addition
$\tau ,\sigma $ have label sequences differing by a 321-avoiding
permutation, then the gradient path from $\tau $ to $\sigma $ is 
unique.
\end{cor}

\proof
A reduced expression for a 321-avoiding permutation 
cannot have consecutive reflections $s_{i+1}s_is_{i+1} $ or 
$s_is_{i+1}s_i$.  This eliminates the possibility of two choices for
the conclusion of the gradient path.
\EOP

\begin{rk}
It would be nice if Theorem ~\ref{red_exp} could be extended
to deal with gradient
paths on nonsaturated chain segments, at least for 
least-content-increasing labellings.  
This is related to the question of describing incidences 
in the complex $\Delta^M$ of critical cells for a 
lexicographic discrete Morse function $M$.
\end{rk}

\medskip
In [13], the preceding theorem is applied to monoid posets to provide
connectivity lower bounds in terms of Gr\"obner basis degree.
Monoid posets are not Cohen-Macaulay in general, so something other
than shelling was needed.  

\begin{examp}\label{grad-examp}
\rm{
Recall the poset from Example ~\ref{lex-examp}.
Consider the saturated chains $F_j$ and $F_i$, labelled $26543$ and 
$23654$, respectively.  $F_j$ has minimal skipped intervals with rank
sets $\{ 1\} , \{ 2\} , \{ 3\} , \{ 4\}$.  The first interval results from a 
syzygy, and the others come from descents.  The facet $F_i$ was
already discussed in Example ~\ref{lex-examp}.
$F_j$ and $F_i$ contribute critical cells
$\tau = x_2 < x_2x_6 < x_2x_6x_5 < x_2x_6x_5x_4$ and 
$\sigma = x_2 < x_2x_3x_6 < x_2x_3x_6x_5$, respectively.  
There is a gradient path $d_4\circ u_4 \circ d_3 \circ u_3 \circ d_2 
= s_4 \circ s_3 \circ d_2$ from $\tau $ to $\sigma $, shown in Figure
~\ref{fig-gradpath}.
\begin{figure}[h]
\begin{picture}(250,160)(0,0)
\includegraphics[width=0.7\textwidth, angle = 0]{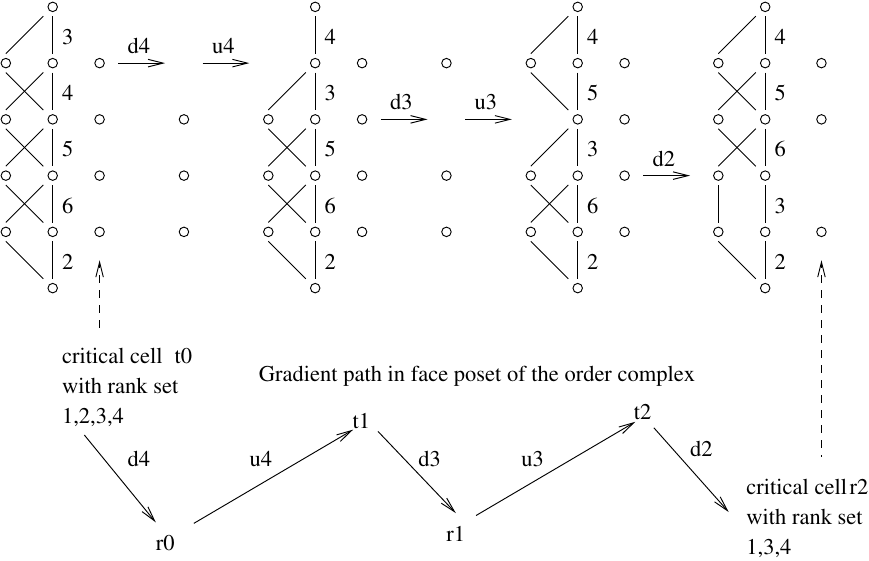}
\end{picture}
\caption{A gradient path between two critical cells}
\label{fig-gradpath}
\end{figure}

It is not hard to check the 
conditions of Theorem ~\ref{red_exp}, using the fact that for any 
elements $u_1\prec u_2 \prec u_3 $ with descending labels 
$\lambda (u_1,u_2) > \lambda (u_2,u_3)$ from the label set $\{ x_3,\dots 
x_6\} $, there exists 
$u_1\prec z \prec u_3$ with $\lambda (u_1,z) = \lambda (u_2,u_3)$ and 
$\lambda (z,u_3) = \lambda (u_1,u_2)$.
Theorem ~\ref{red_exp} 
ensures that the gradient path constructed above is the only
gradient path from $\tau $ to $\sigma $.
The point is that $x_3$ is {\bf non-essential} to the syzygy in 
$\sigma $, allowing it to be shifted to its location in $\tau $ above
the syzygy interval.  Theorem ~\ref{red_exp} 
ensures that there is a unique way for a gradient path to
shift the label $x_3$ downward from its position in $\tau $ to inside
the syzygy interval in $\sigma $.
}
\end{examp}

Next is an example with two gradient paths from $\tau^{(p+1)}$
to $\sigma^{(p)}$, where we examine 
orientation to determine that $\tau^{(p+1)} $ 
has 0 boundary due to cancellation.

\begin{examp}\label{two-paths}
\rm{
Consider the semigroup ring $k[x_1,\dots ,x_7]/(x_4x_5x_6-x_1x_2x_3) \sim
k[ab,cd,ef,ad,be,cf,g]$ and the partial order on its monomials by
divisibility.  As before, label poset edges with the indices of the 
variables being multiplied.  Let us examine the interval $(1,x_4x_5x_6x_7)$.
It has a critical 1-cell $\sigma $
in the saturated chain labelled $7,4,5,6$ and it has critical 2-cells
$\tau_1,\tau_2$ in the saturated chains labelled $7,3,2,1$ and $7,6,5,4$,
respectively.  The order complex has a 2-sphere coming from $\tau_1$, and
then attaching $\sigma $ gives a wedge of a 2-sphere with a 1-sphere.  
When we next attach $\tau_2$, notice that there are two gradient paths
from $\tau_2$ to $\sigma $, and $\sigma $ is incident to $\tau $
in two different ways (coming from two distinct series of collapses).  To
decide whether attaching $\tau_2$ creates a 2-sphere or a 
copy of $\reals P_2$, we need to examine relationships between 
orientations (see [10] for more about orientation).
The orientations are such that $\tau_2$ has 0 boundary, due
to cancellation, and we deduce that $\Delta (1,x_4x_5x_6x_7)$ consists of
a pair of 2-spheres, joined at two distinct points, and has Betti
numbers $\tilde{b}_0=0,\tilde{b}_1=1,\tilde{b}_2=2$.  
}
\end{examp}  

\begin{prop}
If an edge labelling is 
least-content-increasing, then Theorem ~\ref{red_exp} 
holds for $\tau ,\sigma $ of equal content without needing to assume the 
ordered zigzag property or that $\omega $ acts effectively.
\end{prop} 

\proof
Since the initial and final saturated chains have the same label
content, and content is non-increasing, it 
must be preserved at each step.  This together with the 
fact that the labelling is least-increasing implies that each
gradient path downward step which changes the 
label sequence must sort labels on the interval where the
chain element was deleted.  Hence, $\omega $ acts effectively.

Assuming $\tau $ and $\sigma $ agree above rank $s$
means in particular
both have the same set of labels in the same order above rank 
$s$.  Any gradient path step involving ranks above $s$ must sort labels
from this fixed order above rank $s$, precluding completion of the gradient
path to $\sigma $, since gradient path steps can never create new inversions
above rank $s$.
Thus, we are assured that all ranks above $s$ are left untouched by all
gradient paths from $\tau $ to $\sigma $.
\EOP

\begin{figure}[h]
\begin{picture}(250,360)(-15,0)
\includegraphics[width=0.6\textwidth, angle = 0]{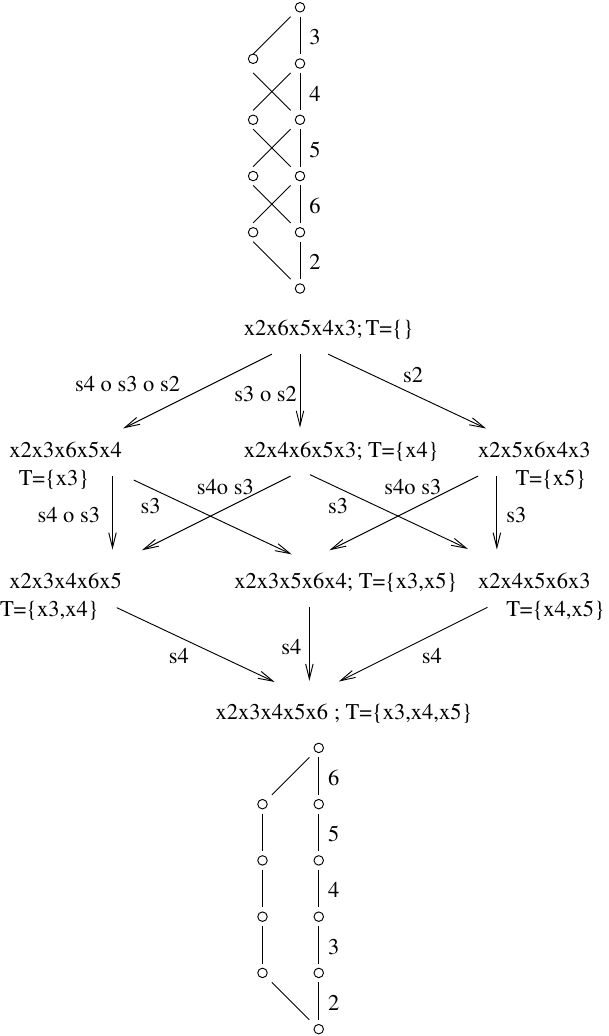}
\end{picture}
\caption{A Boolean algebra of critical cells in $P^M$}
\label{fig-boolcrit}
\end{figure}

Recall from Example ~\ref{lex-examp} the saturated chain 
$F_j$ labelled $x_2x_3x_6x_5x_4$, contributing the
critical cell $\sigma = x_2 < x_2x_3x_6 < x_2x_3x_6x_5$.
In Example ~\ref{grad-examp}, 
we used reduced expressions for permutations (and Theorem ~\ref{red_exp}) 
to show there is a unique gradient path
from $\tau = x_2 < x_2x_6 < x_2x_5x_6 < x_2x_4x_5x_6$ to $\sigma $, 
given by the reduced expression $s_4\circ s_3\circ s_2$.  Now we 
use Theorem ~\ref{acyclic} 
to reverse several gradient paths simultaneously by 
collecting critical cells into a Boolean algebra within $P^M$.

\begin{examp}
Each $T \subseteq S= \{ x_3,x_4,x_5\} $ 
corresponds to a critical cell $Crit(T)$ belonging to facet $F(T)$ 
as follows.  $F(T)$ has label sequence $x_2 x_T x_6 x_{(S\setminus
T)^{rev}} $, where $x_T$ lists members of $T$ in increasing
order, and $x_{(S\setminus T)^{rev} }$ lists members of $S\setminus
T$ in decreasing order.  $F(T)$ has one minimal skipped
interval labelled $x_2x_Tx_6$ coming from the syzygy $x_2x_6=x_1^2$.
Each rank above this interval is covered by a minimal skipped interval
of height one coming from a descent.  Theorem 
~\ref{red_exp} implies that the set of critical 
cells $\{ Crit(T)|T\subseteq S\} $ sits 
inside the multi-graph face poset $P^M$ as a Boolean algebra, depicted in
Figure ~\ref{fig-boolcrit}.  
Its covering relations are $Crit(T\cup \{ x_i\} )\prec 
Crit(T)$ for each $T\subseteq S$ and each $ x_i\in S\setminus T$.  
Matching each $T\setminus \{ x_3\} $ with $T\cup \{ x_3\} $ 
gives a complete acyclic matching on the Boolean algebra.  
\end{examp}

\section{The Cohen-Macaulay property for the poset 
$PD(1^n,q)$ of partial 
decompositions of an $n$-dimensional vector space over $F_q$
into 1-dimensional subspaces}\label{pd1nq}

Let $V$ be a finite vector space, i.e. an $n$ dimensional vector space
over a finite field $\field_q$.  
Let $PD(1^n,q)$ denote the poset of partial 
decompositions of $V$ into 1-dimensional subspaces.  That is, the elements
of $PD(1^n,q)$ are collections $l_1,\dots ,l_r$ of linearly independent
lines in $V$, with a maximal element $\hat{1}$ adjoined.  
\begin{defn}
A set of lines is {\bf independent}
if they are spanned by vectors that are linearly independent.
\end{defn}
$PD(1^n,q)$ has covering relations $u\prec v$ if $v$ is obtained from 
$u$ by adding another linearly independent line to the collection.
$PD(1^n,q)$ is the poset of independent sets for the matroid whose 
ground set consists of the 
lines in $V$, i.e. the points in the projective space
$\proj V$.  Thus, $PD(1^n,q)$ is shellable, by virtue of its order 
complex being the barycentric subdivision of a matroid complex (see [5]).
However, Cohen-Macaulayness may also easily be proven 
using discrete Morse theory in a way that is indicative of how our machinery
also seems to apply to much more complex posets.

We will construct a discrete Morse function in terms of an ordering for the 
ground set of lines in $V$, then use gradient path 
reversal to optimize it into one which has
only top-dimensional critical cells.  The surviving critical cells will
be indexed by 
matroid bases with internal activity 0, allowing us to recover the fact
that a matroid complex has the homotopy type of a 
wedge of spheres of top dimension, where the number of spheres is the 
number of matroid bases with internal activity 0.  Though our proof is
for a specific matroid complex, the same argument will work 
without modification for any isthmus-free matroid.


\subsection{Construction of lexicographic discrete Morse function}

First we need an edge-labelling, then we will study the resulting
lexicographic discrete Morse function.  Following matroid theory,
we call the set of all lines in $V$ our {\bf ground set}, and 
choose an ordering $\omega $ on the ground set.  Label each covering
relation $\{ l_1,\dots ,l_r\}  < \{ l_1,\dots ,l_{r+1}\}$ 
with the line $l_{r+1}$ being added to the set of 
linearly independent lines, and use the order $\omega $ on the ground
set to order our labels.  Label covering relations $\{l_1,\dots ,l_n \} 
< \hat{1} $ with a fixed label that is set to be larger than all lines in 
$PD(1^n,q)$.  The resulting lexicographic order on the
label sequences on saturated chains gives rise to the lexicographic
discrete Morse function we will use.

\begin{prop}
This labelling is least-content-increasing.
\end{prop}

\proof
Any
descent comes from a consecutive pair of lines $l_j,l_i$ with the 
larger one (in terms of the order $\omega $)
inserted first.  This insertion order
may be reversed to obtain an earlier saturated chain on the interval
which has the same content.
\EOP

\begin{defn}
Any independent set of lines that spans $V$ is a {\bf basis}.
An element $l$ in a basis $B$ is {\bf internally
active} in $B$ if $l$ cannot be replaced by a smaller element of 
$F$ (with respect to $\omega $) to obtain an alternate basis for 
$F$.  In other words, $l$ is internally active when
there are no smaller elements whose unique
expression in terms of $B$ has nonzero coefficient for $l$.  The
{\bf internal activity} of a basis is the number of internally 
active elements in it.  (See [5] for more detail.)
\end{defn}

Our first task is to characterize minimal skipped intervals.
\begin{defn}
Let $N$ be a saturated chain $\hat{0} \prec x_1 \prec \cdots \prec
x_k \prec x \prec y_1\prec\cdots \prec y_l \prec \hat{1} $ with label
set $B$.  Suppose that $N$
has strictly increasing labels from $x$ to $\hat{1}$.  Then 
$(x,\hat{1})$ is called a {\bf top interval} of $N$ if the following
conditions are both met:
\begin{enumerate}
\item
The label for $ x\prec y_1$ is not internally active in $B$.
\item
The label on each $y_i \prec y_{i+1}$ is internally active in $B$.
\end{enumerate}
\end{defn}

\begin{prop}
The lexicographic discrete Morse function given by our ground set
labelling/ordering has exactly the following
minimal skipped intervals: (1)
label sequence descents, and (2) top intervals. 
\end{prop}

\proof
An interval may only have one increasing chain unless
the interval is of the form $(x,\hat{1})$, since our labelling 
restricted to $(x,y)$ for $y\ne \hat{1}$ is the standard EL-labelling 
on a Boolean algebra.  Thus, all minimal skipped intervals come from
either descents or intervals $(x,\hat{1})$, since our labelling is
least-content-increasing.  

For a 
skipped interval $(x,\hat{1})$ to be a minimal skipped interval in a 
saturated chain $M = \hat{0} \prec m_1\prec m_2\prec \cdots \prec
x \prec y \prec \cdots \prec \hat{1}$, 
there must be an earlier saturated chain $M'$ that 
agrees with $M$ through $x$ but has covering relation $x\prec y'$ with
strictly earlier label $w'$ than the label $w$ on $x\prec y$.  Let
$w_1,\dots ,w_i$ be the labels on covering relations in
the saturated chain $M$ below $x$.  Then $w'$ is not in the span of 
the vectors $w_1,\dots ,w_i$, since $w_1,\dots ,w_i,w'$ extends to a 
basis.  But for the skipped interval 
$(x,\hat{1})$ to be minimal, $w'$ must be in $\langle w_1,\dots ,w_i,w
\rangle $, because otherwise the expression for $w'$ in terms of $B$
would involve some later label, causing that label to not be 
internally active.  This would mean there would be a smaller minimal
skipped interval beginning later in that saturated chain, a contradiction.
Thus, the expression for $w'$ in terms of $B$ must involve $w$, which 
means $w$ may be replaced by $w'$ to obtain a new basis, implying $w$
is not internally active.  If any later label were internally active,
that would contradict the minimality of our minimal skipped interval from
$x$ to $\hat{1}$, so we have shown that any minimal skipped interval not
coming from a descent must be a top interval.
\EOP

\begin{cor}
A saturated chain contributes a critical cell if and only if it has
label sequence $b_mb_{m-1}\cdots b_1ac_1\cdots c_p$ with 
$b_m>b_{m-1}>\cdots >b_1>a$, $a<c_1<\cdots <c_p$, in which $c_1,\dots ,
c_p$ are all internally active but $a$ is not internally active.
\end{cor}

\proof
A saturated chain contributes a critical cell if and only if its
minimal skipped intervals cover all ranks, so the saturated chain
must have descending labels up until a minimal skipped interval
$(x,\hat{1})$.  However, top intervals are the only minimal skipped
intervals of the form $(x,\hat{1})$.
\EOP

\subsection{Cancellation of critical cells}

This section uses a notion called the non-essential set of a saturated
chain to collect critical cells into Boolean algebras within $P^M$.
We will use complete acyclic matchings on these Boolean algebras, whenever
the non-essential set is nonempty, to construct an acyclic matching on
$P^M$.  The key will be to choose the non-essential set in such a way
that any subset $S$ of the 
non-essential set will index a critical cell in such a way that the 
critical cells given by subsets of the non-essential set comprise a 
Boolean algebra within 
$P^M$.  It is important to define the
non-essential set in a way that does not depend on the order in which 
labels appear on a saturated chain, to ensure that all saturated chains
indexed by subsets of the non-essential set will have the same 
non-essential set.

\begin{defn}
The {\bf non-essential set} of a saturated chain $N$, denoted $NE(N)$, 
is the set of internally active labels on $N$. 
\end{defn}

Let $a < b_1 < b_2 < \cdots < b_n$ be the labels on $N$.  Let
$b_S$ denote the label sequence $b_{i_1}b_{i_2}\cdots b_{i_r}$
for $S = \{ i_1,\dots ,i_r\} \subseteq [n]$ with 
$i_1<\cdots <i_r$.  Let $b_{S}^{rev}$
denote the label sequence $b_{i_r}b_{i_{r-1}}\cdots b_{i_1}$.  Let
$S^C$ denote $[n] \setminus S$.
In this notation, each saturated chain that contributes a critical
cell is labelled 
$b_{(S^C)^{rev}} a b_{S}$ for some set of labels and some choice of $S$.

\begin{rk}\label{in-ne-set}
If $N$ contributes a critical cell that is not top-dimensional, then
its top interval includes one or more internally active labels, implying
$NE(N)\ne \emptyset $.
\end{rk}

\begin{prop}
If the saturated chain $N$ labelled $b_{(S^C)^{rev}} ab_S$ contributes a 
critical cell, then so does the saturated chain labelled 
$b_{(T^C)^{rev}}ab_T$ for any $T \subseteq NE(N)$.
\end{prop}

\proof
This follows from our characterization of minimal skipped intervals
along with the fact that permuting the order in which 
1-spaces are created in a saturated chain does
not affect whether they are internally active in the basis $B$.
\EOP

We will refer to the critical cell contributed by the saturated chain
$b_{(T^C)^{rev}}ab_T$ as being indexed by $T$, and denote this as
$Crit(T)$.  
The next lemma will be applied both to $PD(1^n,q)$ and to
the poset discussed in Section ~\ref{pisn}. 

\begin{lem}\label{near-EL}
Let $\lambda $ be a least-content-increasing labelling on a finite
poset with $\hat{0}$ and $\hat{1}$.  Suppose $\lambda $ 
restricted to $(\hat{0},x)$ for each $x\ne \hat{1}$ is a
CL-labelling.  Let $N$ be a saturated chain contributing a critical 
cell.  Let $S$ be the set of labels on $N$, excluding the minimal label,
denoted $a$, appearing on $N$.
Let $T$ be a subset of $S$ such that there is a saturated chain 
labelled $(U^C)^{rev}aU$ which contributes a critical cell $Crit(U)$
for every $U\subseteq T$.  If there is a gradient path from 
$Crit(U)$ to $Crit(U)\cup \{ i\}$ 
for each $U\subseteq T$, $i\in
T\setminus U$, then the critical cells indexed by subsets of $T$ form
a Boolean algebra within $P^M$.
\end{lem}

\proof
Construct a gradient path from $Crit(U)$ to 
$Crit (U\cup \{i \} )$ for any $U\subseteq T$, $i\in T\setminus
U$ by shifting the label $i$ upward from its position in the descending
part of a saturated chain to the interior of the top interval.  That is,
downward steps in the gradient path delete the chain element immediately
above $i$, then upward steps replace the newly uncovered descent by the 
ascent which has the pair of labels sorted.  
To deduce uniqueness of each such gradient path, we 
use the version of Theorem ~\ref{red_exp} given in 
Remark ~\ref{dual_red_exp}, since the labelling on $PD(1^n,q)$ is a 
least-content-increasing labelling.  The fact that there are no other
gradient paths between critical cells of consecutive dimension follows
from the fact that there cannot be a gradient path from 
$Crit (U_1)$ to $Crit (U_2)$ if the saturated chain contributing the 
latter has any inversions not present in the former; if there is any
$j\in U_1\setminus U_2$, then $Crit(U_2)$ will have an inversion between 
$j$ and the minimal label on the saturated chain, and this inversion
will be missing from $Crit(U_1)$.  Thus, we are done.
\EOP

\begin{cor}
The critical cells indexed by subsets of $NE(N)$ form a Boolean
algebra within $P^M$.
\end{cor}

\proof
Lemma ~\ref{near-EL} applies because
$PD(1^n,q) - \hat{1}$ is a simplicial poset with each interval
getting the standard EL-labelling for a Boolean algebra.
\EOP

\begin{prop}
The union of the complete acyclic matchings on these Boolean algebras
is an acyclic matching on $P^M$.
\end{prop}

\proof
The point is to preclude directed cycles involving multiple Boolean
algebras, by showing that once we pass from one Boolean algebra to another,
we may never return to the original Boolean algebra.
The idea is to partially order the Boolean algebras resulting from 
various label contents and notice that any gradient path which passes
from one Boolean algebra to another must be a step downward in our
partial order on Boolean algebras, making it impossible for a cycle
to return to the original Boolean algebra.  

The partial order is on choices of basis for $V$, where by convention
each basis vector must have a leading 1.  One basis $B_1$ will be smaller
than another basis $B_2$ in our partial order if $B_1$ is obtained from
$B_2$ by replacing a single vector by a lexicographically smaller vector.
Notice that a downward step in the face poset on the order complex 
$\Delta (PD(1^n,q))$ deletes a single chain element, either preserving
the earliest coatom possible in an extension to a saturated chain, or
else causing this coatom to decrease in our partial order on bases.
\EOP

\begin{thm}
$PD(1^n,q)$ is homotopically Cohen-Macaulay.  It has the homotopy type
of a wedge of spheres of top dimension, with the number of spheres 
equalling the number of bases for $\field_q^n$ with internal activity 0.
\end{thm}

\proof
It suffices to show that each interval has critical cells only of 
top dimension.  This is done by cancelling critical cells in the
lexicographic discrete Morse function resulting from our edge-labelling.
Theorem ~\ref{acyclic} ensures that we can cancel the
critical cells belonging to our union of complete acyclic matchings on
Boolean algebras, since this is an acyclic matching on $P^M$.
Remark ~\ref{in-ne-set}
ensures that every critical cell that is not top-dimensional belongs to
such a Boolean algebra $B_n$ for $n\ge 1$, and hence to a Boolean algebra
in which all cells are indeed matched and cancelled.
Thus, the interval $(\hat{0},\hat{1})$ is homotopy equivalent to a wedge
of spheres of top dimension.  Note that the surviving critical cells
are labelled by the decreasing chains with empty non-essential set, i.e.
where none of the labels are internally active.  There is one such 
saturated chain for each basis with internal activity 0, yielding the 
homotopy type result.  

The same labelling and argument works for 
intervals $(x,\hat{1})$.  Let $x = l_1,\dots ,l_r$, and note that
the labels given by bases for $l_1,\dots ,l_r$ cannot belong to the 
non-essential set of a saturated chain on this interval.  Intervals
$(\hat{0},y)$ for $y\ne \hat{1}$ are Boolean algebras, and hence are
shellable, so Cohen-Macaulayness follows.
\EOP

\section{The Cohen-Macaulay property for
the poset $\Pi_{S_n}$ of partitions into 
permutation cycles}\label{pisn}

Jeff Remmel recently defined the following permutation-analogue of the
partition lattice (personal communication).
Let $\Pi_{S_n}$ be the poset of permutations in $S_n$, partially
ordered as follows,
with a maximal element $\hat{1}$ adjoined.  A permutation $\sigma $ 
covers a permutation $\tau $ if $\tau $ may be obtained from $\sigma $
by splitting a single cycle of $\sigma $ (written in cycle notation)
into two smaller cycles by the following procedure: let
$a_1,\dots ,a_r$ be the elements of the cycle in $\sigma $
listed in such a way
that $\sigma (a_i) = a_{i+1}$  for each $i$.  Choose some 
$S\subseteq [r]$, and let $\tau $ have cycles $(a_{i_1} ,\dots ,
a_{i_k} )$ with $S = \{ i_1,\dots ,i_k \}$ and $(a_{j_1},\dots ,
a_{j_{r-k}})$ with $[r]\setminus S = \{j_1,\dots ,j_{r-k}\} $.
In other words, $\sigma $ is obtained from $\tau $ by shuffling 
together two cycles of $\tau $ to obtain one larger cycle, with 
any cyclic shift of each cycle allowed prior to shuffling.  
(Shuffling takes place on two labelled
circles, in which shuffling amounts to some interspersing of labels.)

\begin{rk}
The number of 
poset elements at corank $k$ is the signless Stirling number of 
the first kind $c(n,k)$.  
The coatoms are the $(n-1)! $ possible
permutations which are $n$-cycles.
\end{rk}

Remmel previously showed that all intervals $(x,y)$ 
for $y\ne \hat{1}$ are isomorphic to intervals in the partition lattice,
and hence are EL-shellable (personal communication).  
In particular, intervals $(\hat{0},y)$ for
$y$ a coatom have $(n-1)!$ decreasing chains, so have M\"obius function
$\mu (\hat{0},y) = (-1)^{n-1}(n-1)! $.  We will show that $\Pi_{S_n}$ 
is Cohen-Macaulay, but it remains open whether or not it is shellable.
 
\subsection{Construction of lexicographic discrete Morse function}
 
Label covering relations with ordered pairs $(i,\pi )$ where
$i$ is an integer and $\pi$ is a permutation in $S_{[2,n]}$, obtained
as follows.  Each covering relation
$u < v$ for $v\ne \hat{1}$ merges two
cycles $C_1, C_2$ by shuffling them.  
Let $i = \max (\min C_1, \min C_2)$ and let
$\pi $  be the lexicographically smallest
permutation obtainable from $v$ by shuffling the
cycles of $v$ (in the sense described above).
 
Order integers $i$ linearly and order permutations $\pi $
by the lexicographic order on their expressions in
one-line notation.  The integer label takes precedence
over the permutation for ordering labels.  
Label all covering relations 
$u \prec \hat{1}$ with some fixed label that is chosen to be 
larger than all other labels being used.
 
\begin{rk}
Ascents and descents are determined by the first
coordinate, since just like in the partition lattice
saturated chains have the distinct integers $2,\dots ,n$ as the 
first coordinates of their labels.
\end{rk}

\begin{prop}
The labelling is
least-content-increasing.
\end{prop}

\proof
This is immediate for intervals $(x,y)$ with $y\ne \hat{1}$, since 
then we have an interval in the partition lattice, and the first
coordinate gives an EL-labelling in which all label sequences have
equal content.  If $y=\hat{1}$, then the saturated chains with 
smallest content are those whose coatom is the lexicographically 
smallest permutation $\pi $ obtainable by shuffling together the cycles
in $x$.  Each such saturated chain has label set 
$\{ (i,\pi )| 2\le i\le n\}$, i.e. with second coordinate fixed as
$\pi $.  The lexicographically smallest saturated chain on our 
interval is the unique increasing chain with this content, implying 
the result.
\EOP

Next we characterize the minimal skipped intervals.  Recall that
a saturated chain contributes a critical cell
if and only if its interval system covers all ranks,
so we also characterize which saturated chains contribute
critical cells.  By convention, 
list the elements of a cycle in the order in which they eventually 
appear in a coatom of the saturated chain being considered.  

\begin{defn}
Let $N$ be a saturated chain $\hat{0} \prec x_1\prec\cdots\prec x_k
\prec x\prec y_1\prec\cdots\prec y_l \prec \hat{1}$ with 
increasing labels from $x$ to $\hat{1}$.
A {\bf top interval} in $N$ is an interval $(x,\hat{1})$ 
such that (1) there is a coatom above $x$ that is smaller than $y_l$,
and (2) $y_l$ is the smallest coatom above $y_1$
\end{defn}

\begin{prop}
A saturated chain $N$ has two types of
minimal skipped intervals, descents and top intervals.
\end{prop}

\proof
Since the labelling is least-content-increasing,
each descent gives a minimal skipped interval of height one, and 
all other minimal skipped intervals come from delinquent chains.
However, our labelling restricted to any interval $(x,y)$ for $y\ne 
\hat{1}$ is essentially the standard EL-labelling for the partition
lattice, so delinquent chains are only possible on intervals 
$(x,\hat{1})$, and then must take the form above.  
\EOP

Thus, any saturated chain contributing a critical cell must have 
descending labels immediately followed by a top interval.

\subsection{Cancelling critical cells}

We will again use 
non-essential sets to collect
critical cells into Boolean algebras, much as we did in our analysis
of $PD(1^n,q)$.  The 
non-essential set of a saturated chain will depend only on its coatom
permutation $\sigma $ (specifically on the inversions in $\sigma $)
and on the following tree structure associated to any saturated chain.

\begin{defn}
A covering relation merging cycles $C_1,C_2$ gives rise to a 
{\bf merge step} denoted $e_{i,j}$, where $i=\min C_1$ and $j=\min C_2$.
Associate a labelled, rooted tree on vertex set $\{ 1,\dots ,n\} $
to a saturated chain by letting 1 be the root and 
including each edge $e(i,j)$ where $e_{i,j}$ is a merge step in the 
saturated chain.  Notice that labels
are decreasing along tree paths to the root, and
there are $(n-1)!$ such trees.   
\end{defn}

The non-essential set of a saturated chain will be a certain subset of
its merge steps.
It will develop that these two pieces of data associated to a saturated
chain, namely its co-atom and its tree structure, will remain constant on
certain Boolean algebras of critical cells, ensuring that all elements
of such a Boolean algebra have the same non-essential set.  
To decide which merge steps belong to the non-essential set, we will 
need a notion of forced and unforced inversions in the co-atom permutation.
The non-essential set of a saturated chain will consist of those 
merge steps which shuffle two cycles in the  
lexicographically smallest possible way.

Associate to each co-atom $\sigma $ the permutation $o(\sigma )$ obtained
by viewing the unique cycle in $\sigma $ as a permutation in one-line
notation; $\sigma $ is cyclically shifted so that 1 appears in the leftmost
position, implying 1 is a fixed point of the permutatin $o(\sigma )$.

\begin{defn}
A merge step $e_{m_1,m_2}$ in a saturated chain with
co-atom $\sigma $ is said to have an {\bf unforced} inversion 
if $o(\sigma )$ has an inversion $(i,j)$ with $i,j$ belonging
to the two distinct cycles, $C_1,C_2$, respectively, 
to be merged by $e_{m_1,m_2}$, but $C_1$ has no inversions $(i,k)$
with $i<j<k$.  Any other inversions between elements of $C_1,C_2$
are said to be {\bf forced}.  
\end{defn}

The following lemma will help us to define non-essential sets as 
consisting of merge steps with no unforced inversions.

\begin{lem}\label{forced-describe}
An interval $(x,\hat{1})$ in a saturated chain $N$ is a top interval
if and only if it consists of merge steps $(1,i_1), (1,i_2),
\dots ,(1,i_l)$ with $i_1<i_2<\cdots < i_l$ such that $e_{(1,i_1)}$ 
has at least one unforced inversion, but 
all subsequent merge steps $e_{(1,i_j)}$ have only forced inversions.
\end{lem}

\proof
Let 
$\underline{a} = (a_1,\dots ,a_k)$ and 
$\underline{b} = (b_1,\dots ,b_l)$ be two permutation cycles,
cyclically shifted to have their smallest elements listed first.
Now view the sequences $a_1,\dots ,a_k$ and
$b_1,\dots ,b_l$ as two tapes to be read from left to right.  To 
merge the sequences in the lexicographically smallest possible way,
we proceed as follows:
\begin{enumerate}
\item
Compare the current positions on the 
$\underline{a}$ and
$\underline{b}$ tapes at each step.  
\item
Choose the smaller 
of the two values currently being viewed.
\item
Append this number to the
shuffled sequence being constructed.
\item
Move forward one position
on the tape from which we chose a number.  
\end{enumerate}
From this viewpoint,
it is clear that we cannot add $b_j >a_i$ to the merged sequence
before putting $a_i$ into the merged sequence unless there is 
some $a_{i'}>b_j$ appearing earlier than $a_i$ on the 
$\underline{a}$-tape, since
we will not choose $b_j$ until the $\underline{a}$-tape is at a value 
$a_{i'}$ that is larger than $b_j$.  
Thus, for $\underline{a}$ 
and $\underline{b}$ to be merged in the lexicographically smallest 
possible way, all inversions between $\underline{a}$ and $\underline{b}$
must be forced inversions.  

In the other direction, if $\underline{a}$ and $\underline{b}$ are merged
in a way that is not lexicographically smallest possible, consider the 
first place where the merged sequence differs from the lexicographically
smallest choice.  That is, consider the first step where we compare an
element $a_i$ on the $\underline{a}$ tape to an element $b_j$ on the 
$\underline{b}$ tape and choose the larger of the two.  Without loss of
generality, say $a_i < b_j$.  For this to be the first deviation from
the lexicographically smallest choice, we are assured that there are no
$(a_i,c)$ inversions on the $\underline{a}$ tape for $a_i < b_j < c$, 
because once such a $c$ is reached on the tape, we are assured that one
tape or the other will have some $d\ge c$ as its current value until either
we reach the end of one tape or the other, or we deviate from the 
lexicographically smallest choice of how to merge the two tapes.   Thus,
we cannot get to the state of having $a_i,b_j$, which are both smaller 
than $c$, as the current states of the two tapes, without deviating from
the lexicographically smallest merged sequence before this, a 
contradiction.
\EOP

We have already associated a tree structure $T(N)$
to each saturated chain $N$.
To define the non-essential set of $N$, we will also need a set partition
on $[2,n]$.  We obtain this 
partition, denoted $\Pi(N)$, from $T(N)$ 
by deleting all tree edges of
the form $(1,i)$, then taking the graph components as
the partition blocks.

\begin{defn}
The {\bf non-essential set} of $N$ consists of those tree edges
$(1,i)$ such that the component of the partition $\Pi (N) $
containing $i$ does
not have any unforced inversions with other components of 
$\Pi (N) $.
\end{defn}
 
\begin{prop}\label{pres-tree}
Each merge step $(1,i)\in NE(N)$ may appear either in the top interval or
in a unique position below the top interval which causes the labels 
below the top interval to be decreasing.
Shifting the merge step $(1,i)$ from one of these two positions to the 
other does not change the tree structure.
\end{prop}

\proof
The first assertion follows from Lemma ~\ref{forced-describe}.
To see that tree structure is preserved,
we show that the block minimums are preserved, since 
tree edges connect block minimums at merge steps. 
Shifting $(1,i)$ upward to within the top
interval or downward to the decreasing segment of labels
can only change block minimums for blocks whose minimum
element is $i$.  However, any merge steps $(i,j)$ for $j>i$ must
occur below the lower of the two positions where $(1,i)$ may 
appear.
\EOP

\begin{cor}\label{lower-cancel}
Let $N$ be a saturated chain which contributes a critical cell.
Then the labels in the interior of the top-interval of $N$ all belong 
to $NE(N)$.  In addition, any 
subset of $NE(N)$ may appear in the interior of the 
top interval.
\end{cor}

\proof
This follows immediately from 
Proposition ~\ref{pres-tree}.
\EOP

\begin{cor}\label{bool-pm}
Critical cells are arranged into Boolean algebras $B_m$ in $P^M$, 
based on their non-essential sets, with $m$ being the size of the
non-essential set.
\end{cor}

\proof
This is immediate from Lemma ~\ref{near-EL} of the previous section.
\EOP

\begin{rk}
Corollaries ~\ref{lower-cancel} and ~\ref{bool-pm}
allow us to collect all non-top-dimensional critical 
cells, along with some top-dimensional
ones, into Boolean algebras $B_m$ with $m\ge 1$, since all 
non-top-dimensional critical cells have nonempty 
non-essential set.
\end{rk} 
 
Theorem ~\ref{acyclic} shows that
any acyclic matching on $P^M$ gives a collection of pairs
of critical cells that may be cancelled simultaneously. 
Each Boolean algebra $B_m$ for $m\ge 1$ has a complete acyclic matching.  
The next result implies that the union of these matchings 
is an acyclic matching on $P^M$.

\begin{prop}
There are no directed cycles visiting multiple Boolean algebras.
\end{prop}

\proof
Similarly to with $PD(1^n,q)$, we 
partially order Boolean algebras of critical cells and
show that gradient paths may only proceed from later Boolean algebras
to strictly earlier ones in this partial order.  
Any gradient path from one Boolean algebra to another
one either alters the coatom in a way that makes it
lexicographically smaller or alters the tree structure by
turning a pair of edges $(i,j), (j,k)$ into the pair 
$(i,j), (i,k)$.  Thus, we partially order coatoms lexicographically
and order tree structures using covering relations 
$(i,j),(i,k) \prec (i,j),(j,k)$.
Upward steps are
impossible in both partial orders, precluding
return to the original Boolean algebra, i.e. precluding a cycle
involving distinct Boolean algebras.  
\EOP

\begin{thm}
$\Pi_{S_n}$ is a homotopically Cohen-Macaulay poset.
\end{thm}

\proof
The optimized discrete Morse function has only top-dimensional
critical cells, ensuring the interval $(\hat{0},\hat{1})$ is 
homotopy-equivalent to a wedge of spheres of top dimension.  The proof
extends to arbitrary intervals as follows.
Remmel observed that all intervals
$(x,y)$ with $y\ne \hat{1}$ are shellable.  Our discrete Morse function
easily generalizes to intervals $(x,\hat{1})$ by excluding from the 
non-essential set of a saturated chain restricted to $(x,\hat{1})$ any
labels appearing below $x$.  We may still define forced and unforced 
inversions in terms of co-atom permutations.
\EOP

\section*{Acknowledgments}

The author thanks G\"unter Ziegler for helpful comments on an earlier 
version of her paper and for suggesting a related project that 
motivated this one.  She also thanks David Arthur,
Robin Forman, Phil Hanlon, Jeff Remmel, John Stembridge and Volkmar Welker 
for helpful discussions.

\end{document}